\documentclass[conf]{new-aiaa}
\usepackage[utf8]{inputenc}

\usepackage{graphicx}
\usepackage{amsmath}
\usepackage[version=4]{mhchem}
\usepackage{siunitx}
\usepackage{longtable,tabularx}
\title{A Comparison of SOR, ADI and Multigrid Methods for Solving Partial Differential Equations}

\author{Mohamed Mohsen Ahmed}
\affil{Kansas State University, Manhattan, Kansas, USA} 

\begin{document}

\maketitle

\begin{abstract}
This article presents several numerical techniques for solving Laplace equation. A numerical FORTRAN solver is developed to solve the 2D laplace equation. The numerical approaches implemented in the solver include Jacobi, Gauss-Siedel, Successive Over Relaxation, Alternating Direct Implicit and Multigrid methods. Detailed comparison between different numerical methods is presented and discussed.  
\end{abstract}

\section{Nomenclature}

{\renewcommand\arraystretch{1.0}
\noindent\begin{longtable*}{@{}l @{\quad=\quad} l@{}}
$k$& iteration number \\
$L_x$  & domain length in x direction \\
$L_y$  & domain length in y direction \\
$m$  & number of points in x direction \\
$n$ &    number of points in y direction \\
$R$ &    Residual \\
$w$ &    relaxation factor \\
$\psi$ & stream function
\end{longtable*}}

\section{Introduction}
Laplace equation is an elliptic partial differential equation that governs several simple physical problems such as irrrotational incompressible fluid flow and steady state heat transfer in solids. The mathematical formula of a two-dimensional Laplace equation is 
\begin{equation}
\frac{\partial^2 \psi}{\partial x^2} + \frac{\partial^2 \psi}{\partial y^2}=0
\end{equation}

\begin{figure}[hbt!]
\centering
\includegraphics[width=.5\textwidth]{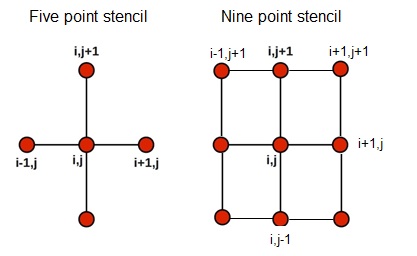} 
\caption{Five and nine point stencils} \label{fig1}
\end{figure}

The finite difference representation of Laplace equation can be expressed in two methods, the five point stencil and the nine point stencil. The geometric interpretation of the two methods is provided in figure ~\ref{fig1}. The five point stencil method is the most common, and its formula is obtained by
\begin{equation}
\frac{\psi_{i-1,j}-2 \psi_{i,j}+\psi_{i+1,j}}{\Delta x^2} + \frac{\psi_{i,j-1}-2 \psi_{i,j}+\psi_{i,j+1}}{\Delta x^2} = 0
\end{equation}
This finite difference formula has a second order truncation error. On the other hand, the nine point formula is obtained by
\begin{multline}
\psi_{i+1,j+1}+\psi_{i-1,j+1} + \psi_{i+1,j-1}+\psi_{i-1,j-1}-2 \frac{h^2-5k^2}{h^2+k^2}(\psi_{i+1,j}+\psi_{i-1,j})
+2 \frac{5h^2-k^2}{h^2+k^2}(\psi_{i,j+1}+\psi_{i,j-1})-20 \psi_{i,j}=0
\end{multline}

This formula is also second order accurate in both spacial directions. However, if $k=h$, the formula becomes of order six in both spacial directions. We are concerned only with the five point stencil formula in the current study.

The solution of the finite difference form of the Laplace equation is similar to the solution of linear algebraic equations. There exist two methods in order to solve any linear group of algebraic equation. The first method is concerned with direct solution of the system of equation. These methods include Crammer's rule, Gauss elimination and Thomas algorithm of tri-diagonal matrices. Although these methods are very simple, they require huge computational time.

Iterative methods, on the other hand, are more efficient than direct methods in most cases since they require far less computational time. In this study, we present an extensive comparison between different iterative methods that are commonly used to solve two-dimensional Laplace equation on a given computational domain \cite{pletcher2012computational}.

\section{Problem Definition}

\begin{figure}[h]
\centering
\begin{minipage}{.5\textwidth}
  \centering
  \includegraphics[width=.9\linewidth]{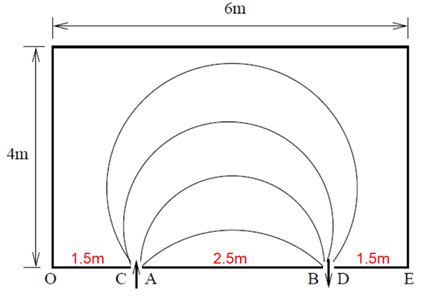}
  \captionof{figure}{Geometry description}
  \label{fig2}
\end{minipage}%
\begin{minipage}{.5\textwidth}
  \centering
  \includegraphics[width=.8\linewidth]{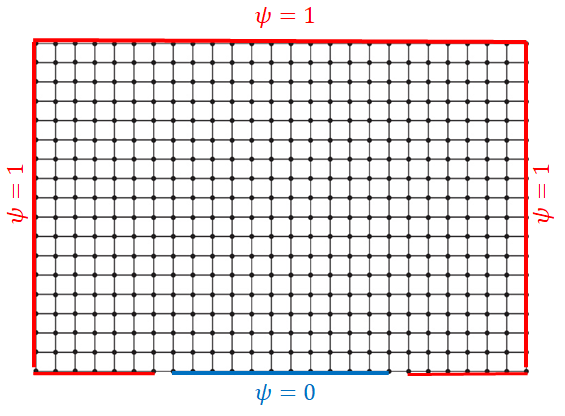}
  \captionof{figure}{Computational grid and boundary conditions}
  \label{fig3}
\end{minipage}
\end{figure}

The problem is concerned with two-dimensional steady incompressible irrational flow in the chamber described in figure ~\ref{fig2}. The parameter $\psi$ here represents the stream function of the potential flow. The inlet and outlet small gaps are $CA=BD=0.25m$. A uniform grid is constructed with constant spacing $\Delta x = \Delta y=0.25 m$. The computational grid and boundary condition is provided in figure ~\ref{fig3}. The stream function $\psi$ has a value of zero on wall AB and zero on all other walls. The mesh size is obtained by
\begin{equation}
m=\frac{L_x}{\Delta x}+1=25 \qquad, \qquad n=\frac{L_y}{\Delta y}+1=17
\end{equation}
where $L_x$ and $L_y$ are the dimensions of the domain in $x$ and $y$ directions respectively. The total number of points in the domain is 425 points including 341 internal points and 84 boundary points.

\section{Iterative Methods}
Iterative methods are used to solve linear algebraic equation by assuming an initial value of the computed parameters and then applying the same algorithm for certain number of calculations in order to converge to the final value of the computed parameters. There exist two types of iterative methods, explicit and implicit. 

\subsection*{Jacobi iterative method}
Jacobi method is an explicit iterative method in which at any iteration cycle every point in the domain is computed from the values of its neighboring points at the old iteration cycle. The finite difference formula employed to compute each point is 
\begin{equation}
\psi_{i,j}^{k+1}=\frac{\psi_{i+1,j}^{k}+\psi_{i-1,j}^{k}+\beta^2(\psi_{i,j+1}^{k}+\psi_{i,j-1}^{k})}{2(1+\beta^2)}
\end{equation}
where $\beta=\frac{\Delta x}{\Delta y}$. Since each element of the right hand side of this equation is at the same iteration cycle (i.e. old iteration), the computations in this method can be parallelized. Therefore, Jacobi method is said to be vectorizable.

\subsection*{Gauss-Seidel method (GS)}
At a current iterative cycle $k+1$ in the calculation of the point $i,j$ using GS method, the neighboring points $i-1$ and $j-1$ that proceed the current point $i,j$ are not considered at the old iteration cycle $k$. Their updated values at the current iteration $k+1$ are used instead of their values at the old iteration $k$.The finite difference formula that is used in GS iterations is expressed as
\begin{equation}
\psi_{i,j}^{k+1}=\frac{\psi_{i+1,j}^{k}+\psi_{i-1,j}^{k+1}+\beta^2(\psi_{i,j+1}^{k}+\psi_{i,j-1}^{k+1})}{2(1+\beta^2)}
\end{equation}

GS iterative method is considered more efficient that Jacobi iterative method in terms of the number of iterations. However, it might not be more efficient in terms of computational time since it is not vectorizable.

\subsection*{Successive over relaxation method (SOR)}

In this method, a relaxation factor $w$ is used in order to accelerate the iterative procedure. SOR method applied to Gauss-Seidel method is

\begin{equation}
\psi_{i,j}^{k+1}=(1-w) \ u_{i,j}^k+ w \frac{\psi_{i+1,j}^{k}+\psi_{i-1,j}^{k+1}+\beta^2(\psi_{i,j+1}^{k}+\psi_{i,j-1}^{k+1})}{2(1+\beta^2)}
\end{equation}
I $1<w<2$ over relaxation is employed and the solution is accelerated. If $w<1$ under-relaxation is employed which makes the solution more stable but slower. If $w>2$ the solution might be unstable. The effect of $w$ on the solution stability and speed is discussed in details in the results and discussion section.

\subsection*{Successive over relaxation by line (SLOR)}
\begin{figure}[h]
\centering
\begin{minipage}{.5\textwidth}
  \centering
  \includegraphics[width=.6\linewidth]{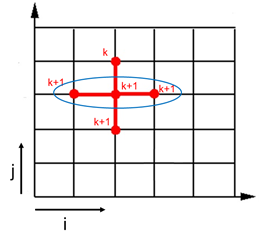}
  \captionof{figure}{SLOR method}
  \label{fig4}
\end{minipage}%
\begin{minipage}{.5\textwidth}
  \centering
  \includegraphics[width=.6\linewidth]{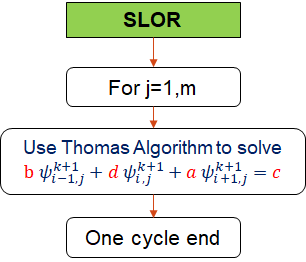}
  \captionof{figure}{Flowchart of single SLOR cycle}
  \label{fig5}
\end{minipage}
\end{figure}

In this method, either rows or column are grouped together as shown in figure ~\ref{fig4}. The solution of each grouped row or column can be obtained implicitly using Thomas Algorithm of tri-diagonal matrices. Figure ~\ref{fig5}  shows the flowchart of a single SLOR cycle.

There are two ways in which relaxation can be employed to update the solution to the next iteration, either before solving by Thomas algorithm (SLORB) or after obtaining the solution from Thomas algorithm (SLORA). The finite difference formula of the SLORA that is solved for each row at $k'+1$ cycle using Thomas algorithm is 
\begin{equation}
-\psi_{i-1,j}^{k'+1}+2(1+\beta^2)\psi_{i,j}^{k'+1}-\psi_{i+1,j}^{k'+1} = \beta^2(\psi_{i,j+1}^{k}+\psi_{i,j-1}^{k+1}
\end{equation}
The solution over relaxation is obtained by
\begin{equation}
\psi_{i,j}^{k+1}=\psi_{i,j}^k+w(\psi_{i,j}^{k'+1}-\psi_{i,j}^k)
\end{equation}

The finite difference equation for SLORB that is used to solve for each row at $k+1$ cycle using Thomas algorithm is

\begin{equation}
-2(1+\beta^2) \psi_{i-1,j}^{k+1}+w \psi_{i,j}^{k+1}-2(1+\beta^2)  \psi_{i+1,j}^{k+1} =2(1+\beta^2)(1-w) \psi_{i,j}^k + w \beta^2(\psi_{i,j+1}^{k}+\psi_{i,j-1}^{k+1})
\end{equation}

\subsection*{Alternative direct implicit method (ADI)}
\begin{figure}[h]
\centering
\includegraphics[width=.255\textwidth]{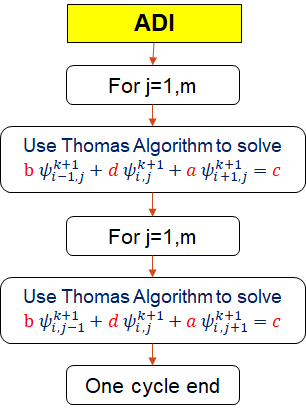} 
\caption{Flowchart of ADI method} \label{fig6}
\end{figure}

ADI method is similar to SLOR method, but in each cycle sweeps by rows are followed by sweeps by columns. Therefore, a complete iteration cycle consists of two steps, as shown in the figure ~\ref{fig6}, which may require more computational time than SLOR method.

\subsection*{Multigrid Method (V-Cycle)}
On fine grids, very large number of iteration may be required to solve Laplace equation using the previously discussed methods. This large number of iteration is mostly cosed by the low frequency component of the error, which can not easily be removed on fine grids. therefore, coarse grids are used in mutligrid method in order to remove the low frequency error. Let the linear operator of the Laplace equation be
\begin{equation}
L(\psi_{i,j})=\frac{\psi_{i-1,j}-2 \psi_{i,j}+\psi_{i+1,j}}{\Delta x^2} + \frac{\psi_{i,j-1}-2 \psi_{i,j}+\psi_{i,j+1}}{\Delta x^2}
\end{equation}
If a solution is obtained using GS method after only three of four iteration (i.e. non converged solution), the right hand side of this equation does not equal zero. And the residual can be obtained by
\begin{equation}
R_{i,j}= L(\psi_{i,j})
\end{equation}
such that $R_{i,j}=0$ at convergence. Let the final converged solution be $\psi_{i,j}$ defined by the correction
\begin{equation}
\psi_{i,j}=\Delta \psi_{i,j} + \psi_{i,j}^k
\end{equation}
where $\Delta \psi_{i,j}$ is the difference between the last two iterations and $\psi_{i,j}^k$ is the solution obtained from iteration $k$. Since $L \psi_{i,j}=0$, we can write
\begin{equation}
L \Delta \psi_{i,j}+ L \psi_{i,j}^k=0
\end{equation}
Using equation 12 we obtain
\begin{equation}
L \Delta \psi_{i,j}=- R_{i,j}
\end{equation} 

\begin{figure}[h]
\centering
\includegraphics[width=.65\textwidth]{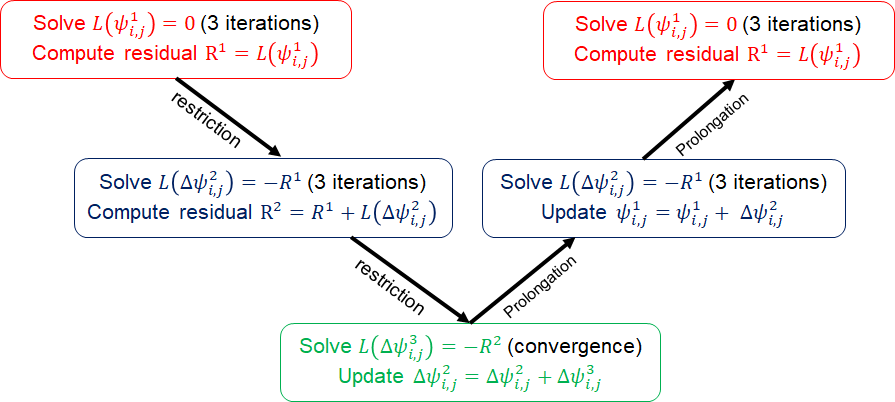} 
\caption{Multigrid V-cycle} \label{fig7}
\end{figure}

Three mesh levels are considered in the multigrid V-cycle. Mesh level 1 has 25x17 points, mesh level 2 has 13x9 points and mesh level 3 has 7x5 points.
The procedure of the V-cycle is shown in ~\ref{fig7}. The solution procedure employed in this study is as follows:
\begin{enumerate}
\item Mesh level 1: solve $L (\psi^1_{i,j})$=0 for 3 iterations by GS method using initial condition $(\psi_{i,j}=0)$
\item Mesh level 1: compute $R^1_{i,j}=L(\psi_{i,j}^1)$
\item Mesh level 2: restrict the values of $R^1_{i,j}$ from mesh level 1 to the coarse grid in mesh level 2
\item Mesh level 2: solve $L(\Delta \psi^2_{i,j})+R^1_{i,j}=0$ for 3 iterations using initial condition ($\Delta \psi^2_{i,j}=0)$
\item Mesh level 2: compute $R^2_{i,j}=R^1_{i,j}+L(\Delta \psi^2_{i,j})$
\item Mesh level 3: restrict the values of $R^2_{i,j}$ from mesh level 2 to the coarse grid in mesh level 3
\item Mesh level 3: solve $L(\Delta \psi^3_{i,j})+R^2_{i,j}=0$ till convergence using initial condition ($\Delta \psi^3_{i,j}=0$).
\item Mesh level 2: interpolate the solution of $\Delta \psi^3_{i,j}$ from mesh level 3 to mesh level 2
\item Mesh level 2: solve $L(\Delta \psi^2_{i,j})+R^1_{i,j}=0$ for 3 iterations using the new initial condition $(\Delta \psi^2_{i,j})_{new}=\Delta \psi^2_{i,j} + \Delta \psi^3_{i,j}$
\item Mesh level 1: interpolate the solution of $\Delta \psi^2_{i,j}$ from mesh level 2 to mesh level 1
\item Mesh level 1: solve $L(\psi^1_{i,j})=0$ for 3 iterations using the new initial condition $(\psi^1_{i,j})_{new}=\psi^1_{i,j} + \Delta \psi^2_{i,j}$
\end{enumerate}
The 11 steps are for one V-cycle. The steps are repeated until the finale solution is converged.

\begin{figure}[h]
\centering
\includegraphics[width=.55\textwidth]{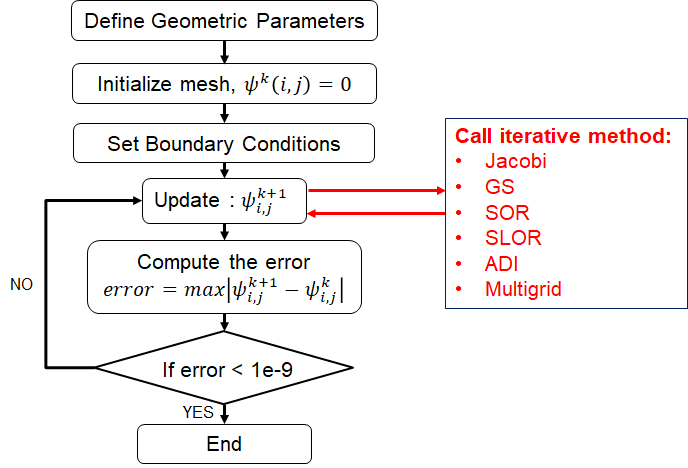} 
\caption{Solver algorithm} \label{fig8}
\end{figure}

A FORTRAN solver is developed to solve the two-dimensional Laplace equation using the six iterative methods discussed previously. A flowchart showing the solver algorithm is described in figure ~\ref{fig8}. The norm infinity error is defined as
\begin{equation}
\textmd{error}=\textrm{max} \mid \psi^{k+1}_{i,j}-\psi^{k}_{i,j} \mid
\end{equation}
The stopping criteria is when the error reaches $10^{-9}$.

\section{Results and Discussion}
\begin{figure}[h]
\centering
\begin{minipage}{.5\textwidth}
  \centering
  \includegraphics[width=.75\linewidth]{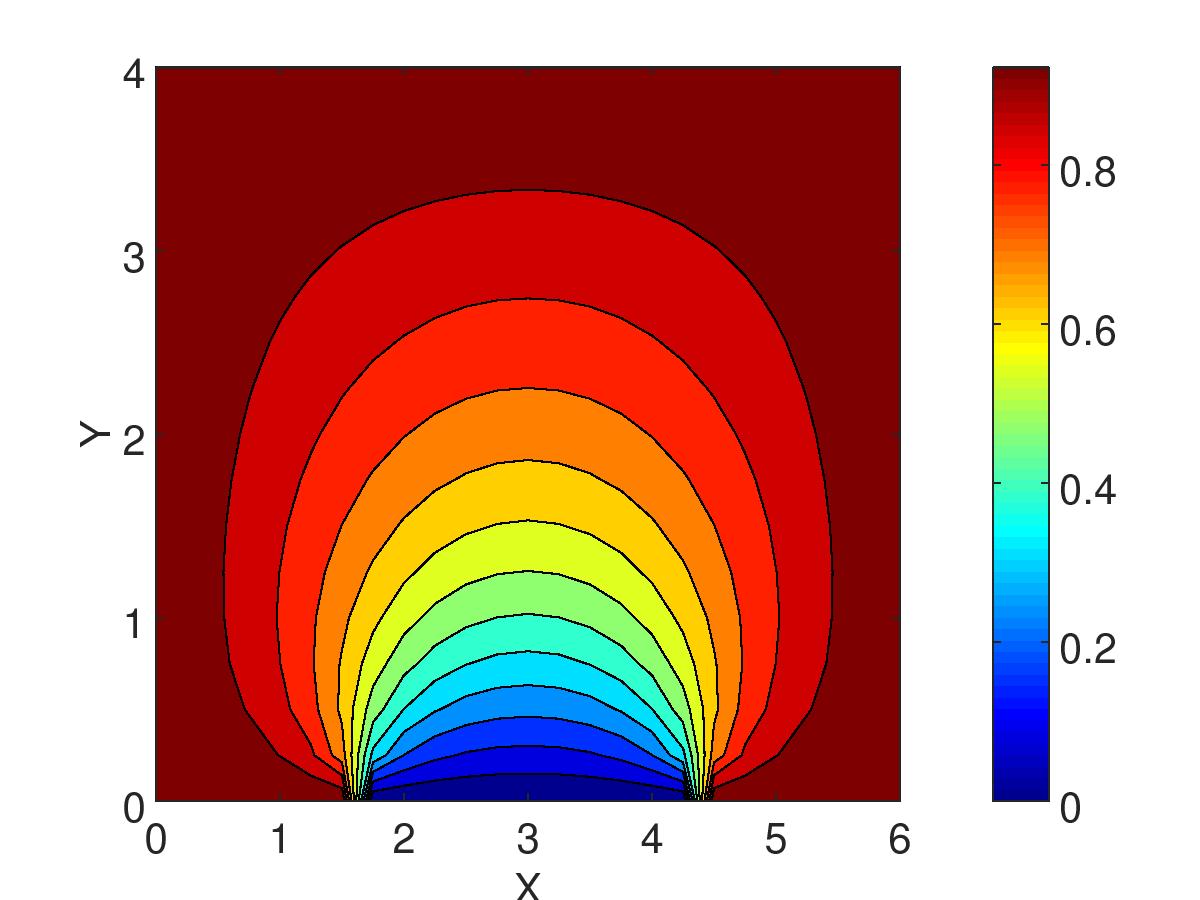}
  \captionof{figure}{Contours of $\psi$ using BCs described in first section}
  \label{fig9}
\end{minipage}%
\begin{minipage}{.5\textwidth}
  \centering
  \includegraphics[width=.75\linewidth]{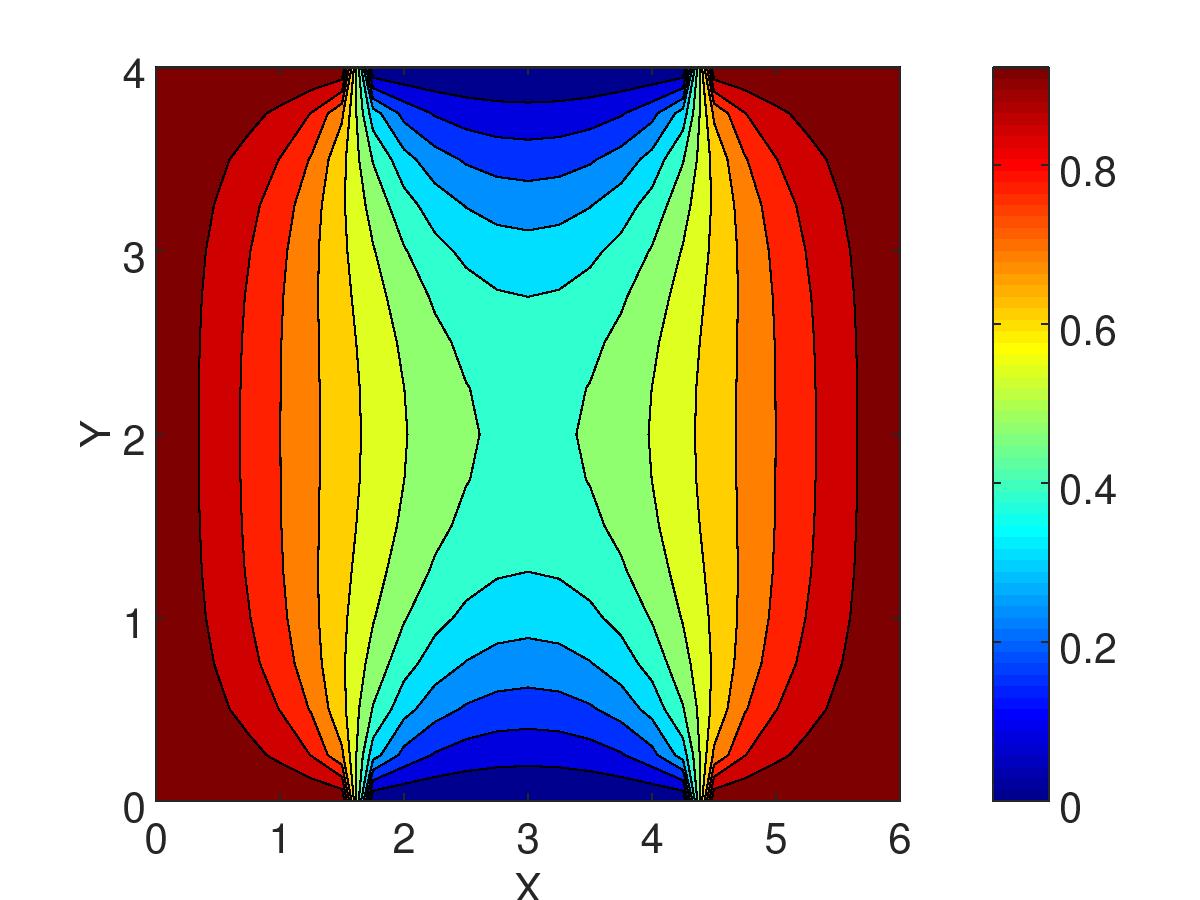}
  \captionof{figure}{Contours of $\psi$ using symmetry Bcs}
  \label{fig10}
\end{minipage}
\end{figure}

The solver is first verified by solving Laplace equation using two different set of boundary conditions. The first set of boundary conditions is similar to that described in the first section. In the second set of boundary conditions, the top wall has similar conditions as the bottom wall. Figures ~\ref{fig9} and ~\ref{fig10} show the contours of $\psi$ obtained from two different boundary conditions. Both solutions give us confidence that the solver is free from coding errors and provide correct physical results.

\begin{figure}[h]
\centering
\begin{minipage}{.5\textwidth}
  \centering
  \includegraphics[width=.9\linewidth]{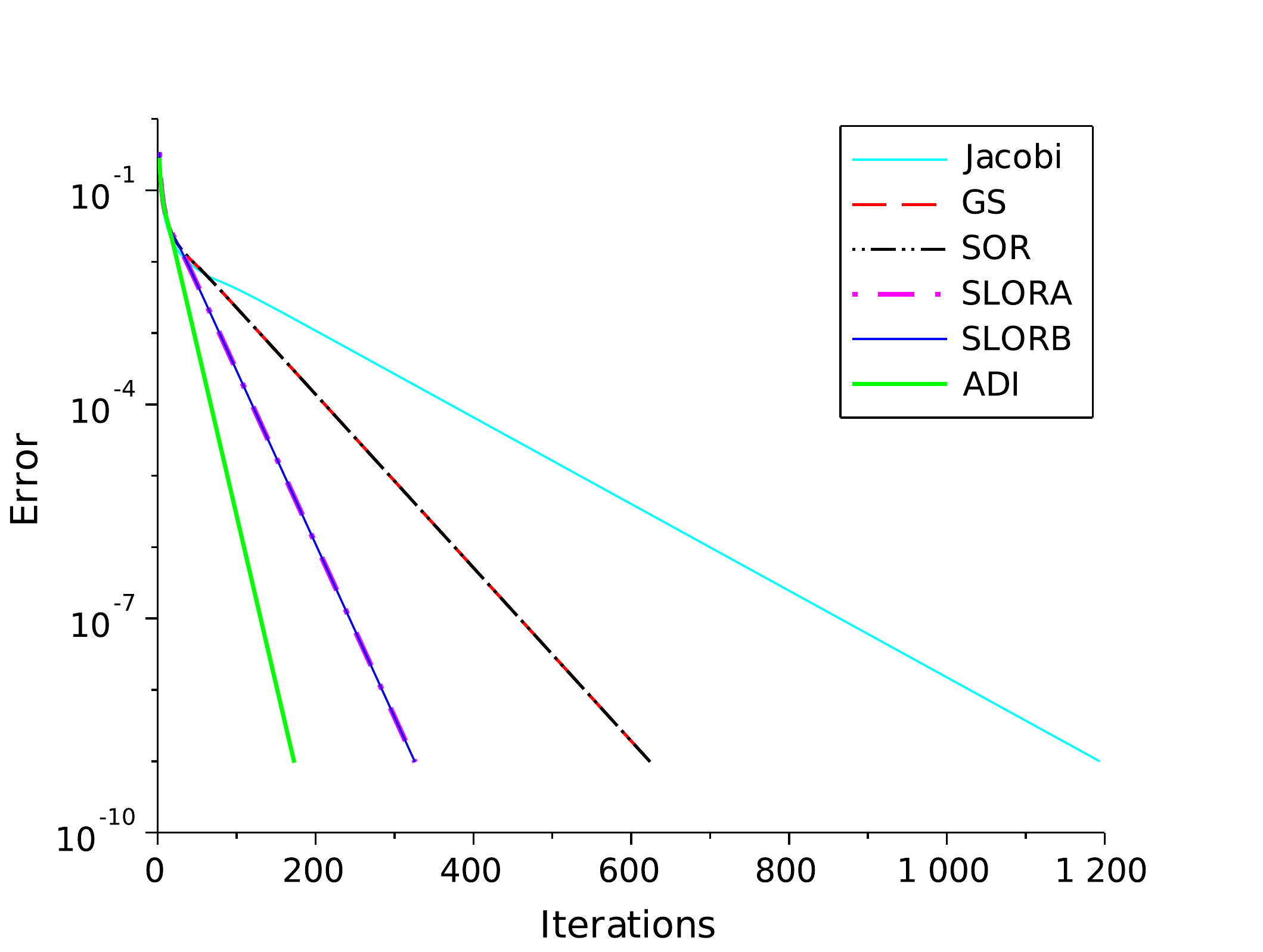}
  \captionof{figure}{Error vs. iterations (no over relaxation)}
  \label{fig11}
\end{minipage}%
\begin{minipage}{.5\textwidth}
  \centering
  \includegraphics[width=1.0\linewidth]{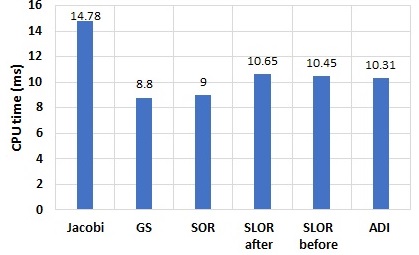}
  \captionof{figure}{Computational time (no over relaxation)}
  \label{fig12}
\end{minipage}
\end{figure}

All numerical methods are first compared without over relaxation in order to give us confidence on the validity of our solver to predict well known trends. Figure ~\ref{fig11} shows the number of iterations and the error of the different iterative methods employed in the solver. The slowest method is Jacobi which requires 1194 iterations to reach convergence. GS and SOR with no relaxation are identical and the number of iterations is 624 iterations, which is approximately half the number of iterations of Jacobi method. SLORB and SLORA are much faster than GS and SOR method with only 326 iterations. Although ADI method is the fastest among these methods with only 326 iterations, it should be taken into consideration that two calculations are conducted inside a single cycle of the ADI method. In figure ~\ref{fig12}, the CPU time for each method is reported. It can be inferred that Jacobi is the slowest method, GS is the fastest method and SLOR and ADI methods are approximately equal in terms of computational time due to the reasons discussed previously. 

\begin{figure}[h]
\centering
\begin{minipage}{.5\textwidth}
  \centering
  \includegraphics[width=1\linewidth]{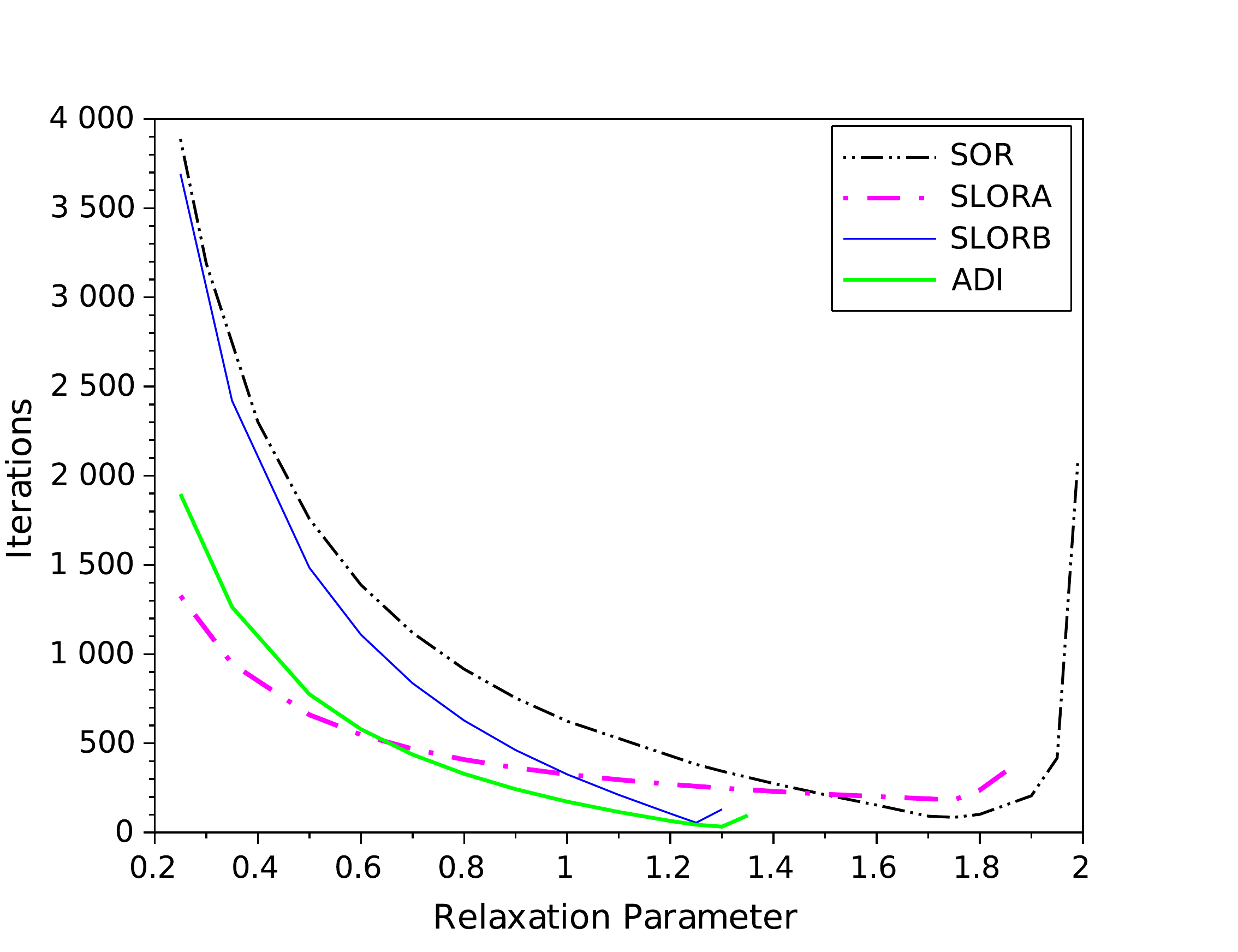}
  \captionof{figure}{iterations vs. relaxation parameter $w$}
  \label{fig13}
\end{minipage}%
\begin{minipage}{.5\textwidth}
  \centering
  \includegraphics[width=1\linewidth]{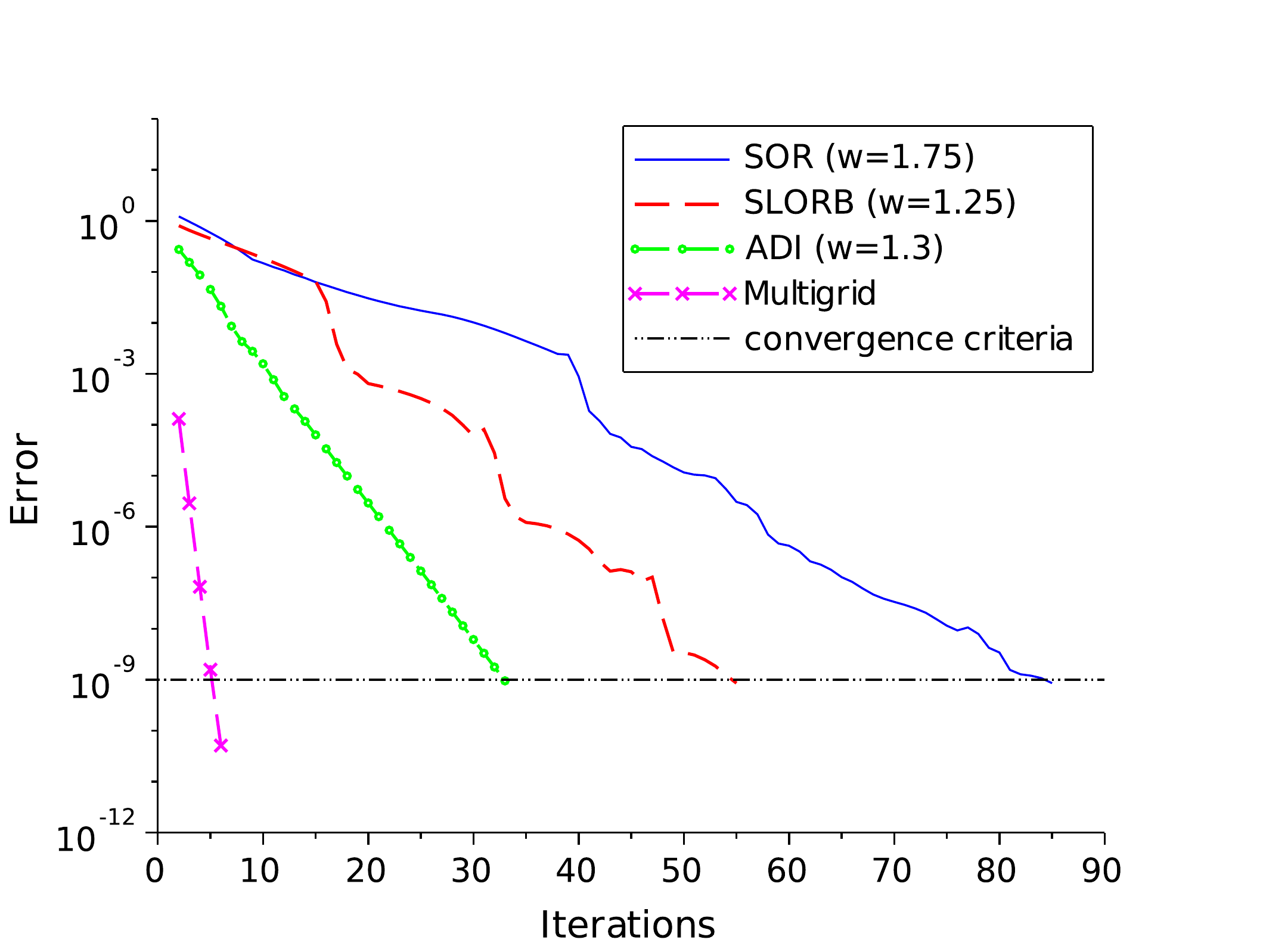}
  \captionof{figure}{Optimum relaxed methods}
  \label{fig14}
\end{minipage}
\end{figure}

Now, we are in turn to find the value of the best over relaxation parameter that would speed up the calculations of the SOR, SLOR and ADI methods. Figure ~\ref{fig13} show the effect of the relaxation parameter on the number of iterations. It can be seen that there exist an optimum value of $1<w<2$ and this value is different for different iteration methods. The optimum values of $w$ are 1.75, 1.75, 1.25 and 1.3 for SOR, SLORA, SLORB and ADI respectively. Figure ~\ref{fig14} shows the error vs. number of iterations required at the best relaxation parameter $w$ for each iterative method. \\

\begin{figure}[h]
\centering
\includegraphics[width=.8\textwidth]{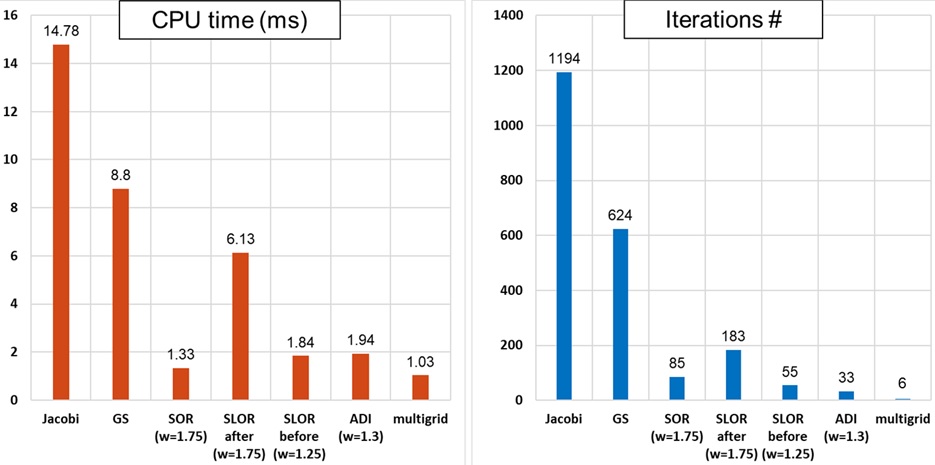} 
\caption{Summary of Iteration number and CPU time} \label{fig15}
\end{figure}

\begin{figure}[h]
\centering
\includegraphics[width=.5\textwidth]{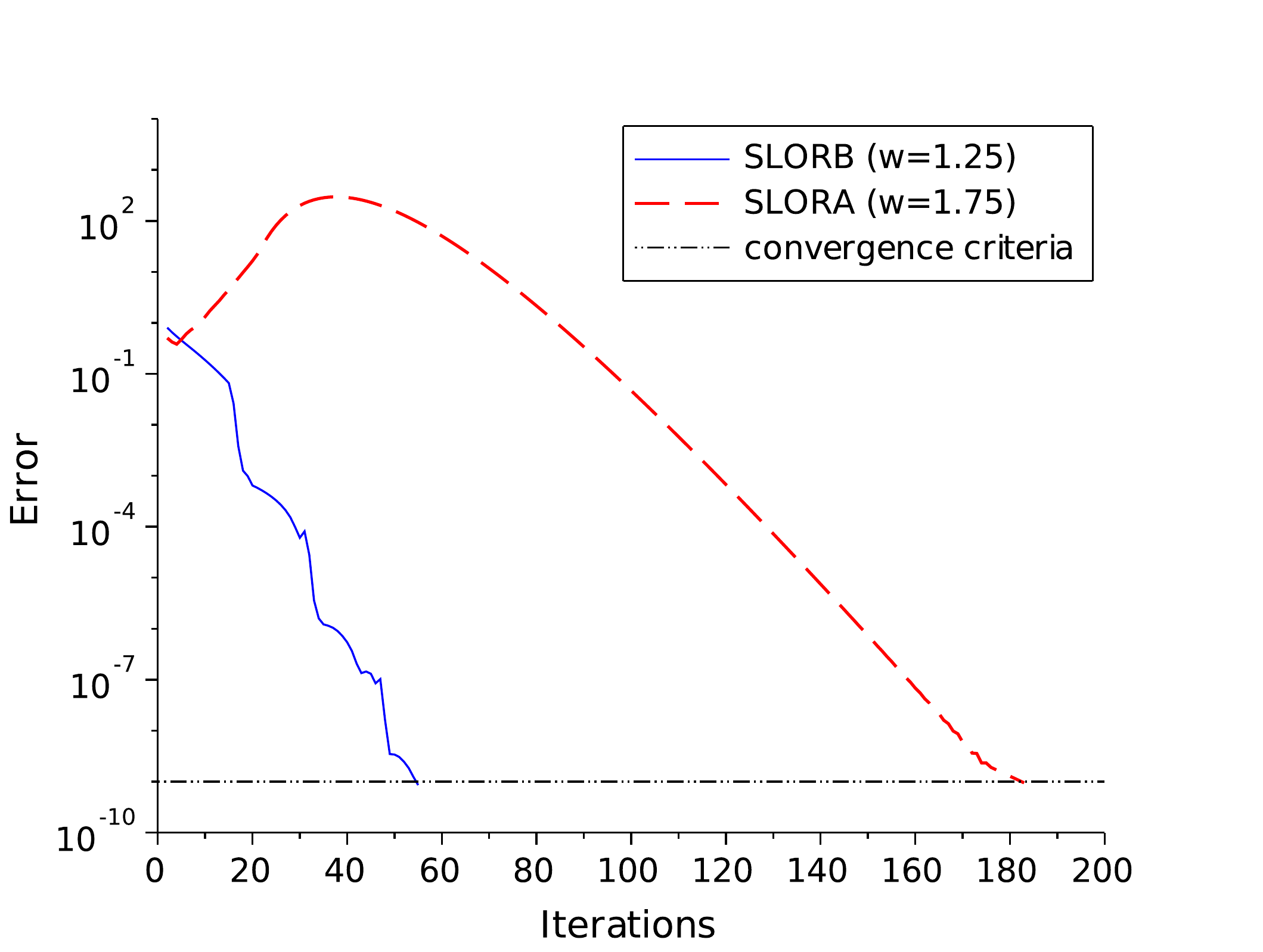} 
\caption{SLORB and SLORA with optimal relaxation parameter} \label{fig16}
\end{figure}

The number of iterations and CPU times for all methods with optimal relaxation are summarized in figure ~\ref{fig15}. As expected, ADI requires less number of iterations that SLORB (33 vs. 55). However, the total CPU time is almost identical in both cases. Another important observation is that SLORA is much slower than SLORB and SOR with optimal over relaxation parameter. Figure ~\ref{fig16} compares SLORB and SLORA at their optimal over relaxation parameters. It is inferred that although $w$ is higher in SLORA, more number of iterations is required since the solution overshoots in the first 50 iterations. The only explanation of this observation is low stability of Thomas algorithm while solving the tri-diagonal matrix when the matrix is not diagonally dominant.

\begin{figure}
\centering
\includegraphics[width=0.95\textwidth]{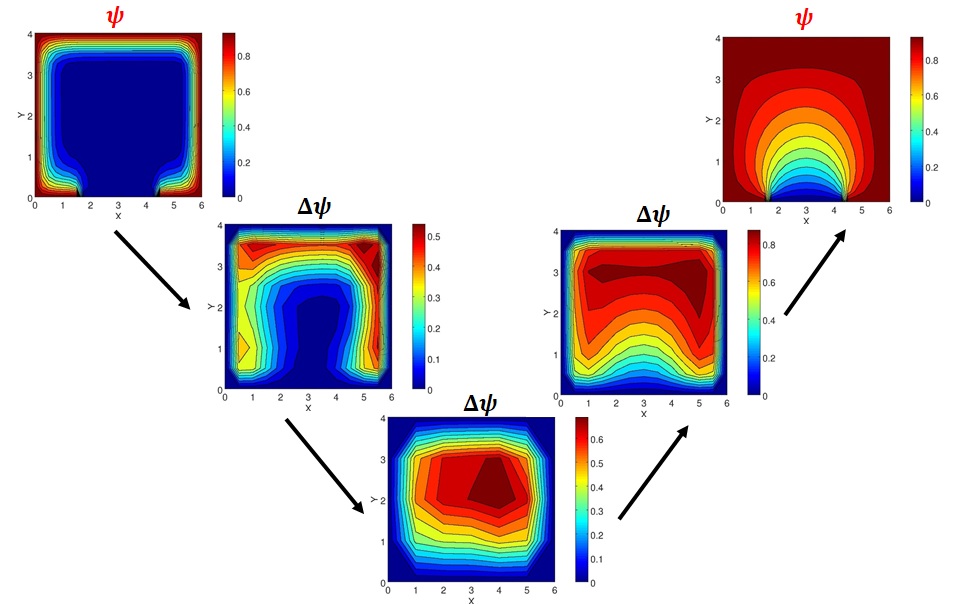} 
\caption{Multigrid V-cycle (cycle 1)} \label{fig17}
\end{figure}

\begin{figure}
\centering
\includegraphics[width=0.95\textwidth]{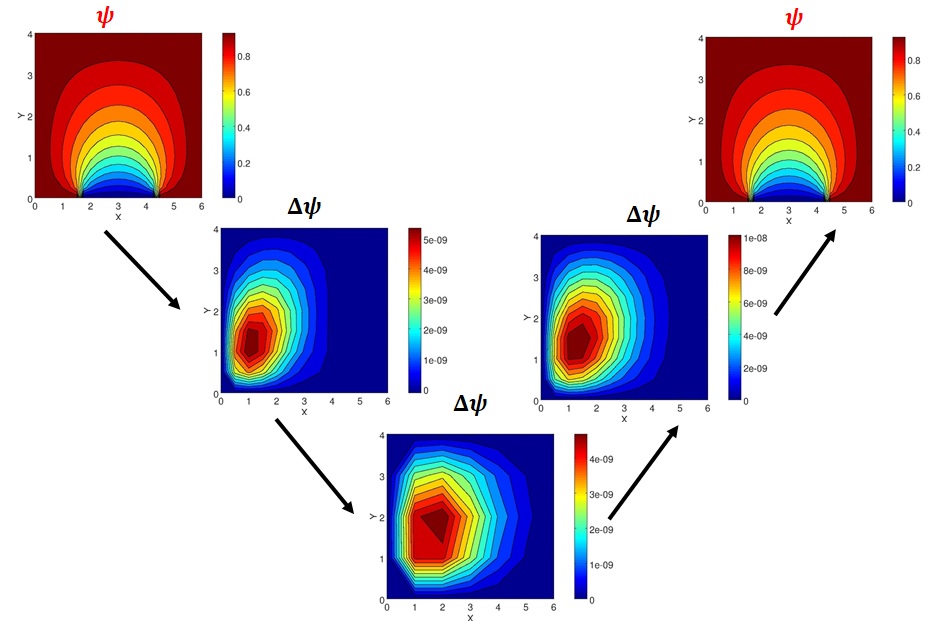} 
\caption{Multigrid V-cycle (cycle 6)} \label{fig18}
\end{figure}

It is also concluded from figure ~\ref{fig15} that multigrid method is the most efficient method among all other iterative methods. The solution converged after only six V-cycles and the total CPU time is around 1.03 milliseconds which is fastest relative to all other methods. The initial and final cycles of the multigid V-cycle method are shown in figures ~\ref{fig17} and ~\ref{fig18} respectively. It is clearly seen in cycle 1 that $\Delta \psi$ is in the order of 1. However, at cycle 6, $\Delta \psi$ is reduced to $10^{-9}$ which is the convergence criteria defined in this study.\\

\section{Conclusion}
Laplace equation is an elliptic second order differentioal equation that describes the behavior of incompressible rotational inviscid flow and heat transfer in solids. Five point stencil method is a second order finite difference scheme of the two-dimensional Laplace equation. The solution of the five stencil finite difference formula is obtained using six different implicit and explicit iterative methods. Over relaxation is employed in order to accelerate the iterative procedure. Comparison between different iterative methods reveal that ADI and SLOR methods are more efficient that other iterative methods in terms of number of iterations and computational time. Multigrid V-cycle method is also employed in the current study and the solution is obtained after only six V-cycles.

\section{Acknowledgment}
The author would like to thank Dr. Mingjun Wei for his guidance and contribution to this work.
\bibliography{sample}

\end{document}